\documentclass[10pt, reqno]{amsart}

\usepackage{amsmath,amsfonts,amsthm,
amssymb,amsxtra,mathrsfs, mathtools, bm, upgreek,
 graphicx,color,caption,subcaption}
\usepackage[pdfborder={0 0 0}]{hyperref}
\usepackage{cases}

 \usepackage[english]{babel}

\newtheorem{Thmk}{\textsc{Theorem}}[section]
\newtheorem{rem}{\textsc{Remark}}[section]
\newtheorem{COROA}{\textsc{Corollary}}[section]
\newtheorem{Prop}{\textsc{Proposition}}[section]
\newtheorem{defi}{\textsc{Definition}}[section]
 
 \let\bb\mathbb
\def\R {{\relax{\bb R}}}

\def\N {\relax{\bb N}}

\newcommand{{\vpz}}{{\varphi^0}}
\newcommand{\vq}{v}
\newcommand{\f}{\frac}
 
\numberwithin{equation}{section}
\def\myAbstract#1{
{\vspace{0.2cm}\flushleft\sffamily\bfseries\Large  Abstract 
\vspace{0.2cm}} 
\\
#1}
 \def\myTalk#1{
\section*{#1}
}
 
 \def\myAuthor#1#2#3#4{
{\noindent #1 #2 \index{#2@#2, #1} \\ #3 \\ #4 \\}}

\def\myCoauthor#1#2#3#4{  
{\noindent #1 #2 \index{#2@#2, #1} \\ #3 \\ #4 \\}}

\def\myLoc#1 {{\flushleft{\bf Chairman:   #1 \\}}}

\begin{document}

\myTalk{\Large \bfseries Existence and stability of a blow-up solution with a new prescribed behavior for a heat equation with a critical nonlinear gradient term}

\bigskip
\myAuthor{Slim}{TAYACHI}
{Universit\'e de Tunis El Manar, Facult\'e des Sciences de Tunis, D\'epartement de
Math\'ematiques, Laboratoire  \'Equations aux D\'eriv\'ees
Partielles LR03ES04,  2092 Tunis, Tunisia}
{slim.tayachi@fst.rnu.tn}

\myCoauthor{Hatem}{ZAAG\footnote{This author is supported by the ERC Advanced Grant no. 291214, BLOWDISOL and by the ANR project ANA\'E ref. ANR-13-BS01-0010-03}}
{Universit\'e Paris
13, Sorbonne Paris cit\'e, Institut Galil\'ee, CNRS UMR 7539 LAGA,
99 Avenue Jean-Baptiste Cl\'ement 93430 Villetaneuse, France}
{Hatem.Zaag@univ-paris13.fr}

\myAbstract{}
 We consider the semilinear heat equation, to which we add a nonlinear gradient term, with a critical power. We construct a solution which blows up in finite time. We
also give a sharp description of its blow-up profile.  The proof  relies on the reduction of the problem to a finite dimensional one,  and uses the index theory to conclude.  
 Thanks to the interpretation of the parameters of the finite-dimensional problem in terms of the blow-up time and point, we also show the stability of the constructed solution with respect to initial data. This note presents the results and the main arguments. For the details, we refere to our paper \cite{TZ15}.
 \\[5mm]
\selectlanguage{english}
\textbf{\textsf{Existence et stabilit\'e d'une solution explosive avec un nouveau comportant prescrit pour une \'equation de la chaleur avec un terme en gradient non lin\'eaire et critique}}.
On consid\`ere l'\'equation semi-lin\'eaire de la chaleur, \`a laquelle on rajoute un terme non lin\'eaire en gradient, avec puissance critique. On montre l'existence d'une solution explosant en temps fini uniquement \`a l'origine, et on en donne le profil \`a l'explosion. Notre m\'ethode s'appuie sur la r\'eduction du probl\`eme en dimension finie, puis la solution de ce probl\`eme gr\^ace \`a un argument topologique. Gr\^ace \`a l'interpr\'etation des param\`etres du probl\`eme de dimension finie en terme du choix du temps et du point d'explosion, on obtient la stabilit\'e de la solution construite par rapport aux donn\'ees initiales. Cette note pr\'esente les r\'esultats et les arguments de la preuve. Pour les d\'etails, voir notre papier \cite{TZ15}.

 \bigskip
 \noindent
{\bfseries AMS 2010 Classification}. {\small Primary
35K55, 35B44; Secondary 35K57.}

\medskip

\noindent {\bfseries Keywords}. {\small Nonlinear heat equation,
blow-up, nonlinear gradient term, blow-up profile, final profile,
single point blow-up, Hamilton-Jacobi equation, asymptotic
behavior  }

\bigskip

\section{Introduction and statement of the results}
\setcounter{equation}{0}

In our paper \cite{TZ15}, we consider the following nonlinear heat equation:
\begin{align}\label{eq:uequation}
\partial_t u & = \Delta u +\mu|\nabla u|^q+|u|^{p-1}u, \\
& u(\cdot,0) = u_0\in W^{1,\infty}(\mathbb{R}^N),\notag
\end{align}
 where $u=u(x,t)\in \R,\; t\in [0,T),\; x\in \R^N$, and the parameters $\mu,\;
 p$ and $q$ are such that
\begin{align}\label{condition:h}
     \mu>0,\; p>3,~  q=q_c\equiv\frac{2p}{p+1}.
     \end{align}
This equation was first introduced by Chipot and Weissler \cite{CWbngsj89} in 1989 with $\mu<0$ (and $q>1$) as a model to understand whether blow-up may be prevented by the addition of the negative gradient term. Later in 1996, Souplet suggested a population dynamics interpretation for the equation in \cite{Smmas96}. Many authors dealt with the mathematical analysis of this equation, both for large time dynamics and finite time blow-up, and also for the elliptic version (see \cite{CWbngsj89,CFQ,Sejde2001,{Snepp2005},SnTW,SnT, ST,STWesbugiumj96,SW} and references therein).

\smallskip
Equation  \eqref{eq:uequation} enjoys two limiting cases:\\
- when $\mu=0$, we have the well-known {\it semilinear heat equation}:
\[
\partial_t u = \Delta u +|u|^{p-1}u;
\]
- when $\mu\to \infty$, we recover (after appropriate rescaling) the {\it diffusive Hamilton-Jacobi equation}:
\[
\partial_t u = \Delta u +|\nabla u|^q.
\]
The value $q=2p/(p+1)$ is critical, since for $\mu \neq 0$, it is the only value for which 
equation \eqref{eq:uequation} is  invariant under the transformation: $u_\lambda(t, x) = \lambda^{2/(p-1)}u(\lambda^2t, \lambda x),$ as
for the equation without the gradient term, that is when $\mu=0$. Moreover, we know from the literature that both the blow-up and the large-time behaviors depend on the position of $q$ with respect to $\frac{2p}{p+1}$ (see the above-mentioned literature).

\smallskip
Equation \eqref{eq:uequation} is wellposed in  $ W^{1,\infty}(\mathbb{R}^N)$ thanks to a fixed point argument (see also \cite{AW92}, \cite{SW} and \cite{SW1}). In our paper \cite{TZ15}, we focus on the study of blow-up for that equation.

\smallskip
When $\mu=0$, there is a huge literature about the subject, and no bibliography can be exhaustive. Let us focus on the existence a stable solution $u(x,t)$ which blows up in finite time $T>0$ only at the origin and satisfies 
\begin{equation}\label{profile0}
(T-t)^{1/(p-1)}u(z\sqrt{(T-t)|\log(T-t)|},t) \sim f_0(z)\mbox{ as }t\to T,
\end{equation}
where 
 {$$f_0(x)=\left(p-1+b_0|x|^2\right)^{-1/(p-1)}\mbox{ and } b_0=(p-1)^2/(4p).$$}
Formal arguments for the existence of such a profile were first suggested by Galaktionov and Posashkov \cite{GP1,GP2} in 1985, then, Berger and Kohn \cite{BergerKhon} gave a numerical confirmation in 1988, and the proof came from Herrero and Vel{\'a}zquez \cite{HViph93} in 1993 and Bricmont and Kupiainen \cite{BK} in 1994. Later, Merle and Zaag \cite{MZsbupdmj97} simplified the proof of \cite{BK} and proved the stability of the profile $f_0$.\\
The authors in \cite {BK} and   \cite{MZsbupdmj97} used a {\it constructive} proof, based on:
\begin{itemize}
\item[-] \vskip-0.2cm
The reduction of the problem to a finite-dimensional one;
\item[-] The solution of the finite-dimensional problem thanks to the degree theory.
\end{itemize}
Let us add that other profiles are possible (see \cite{BK} and \cite{HViph93}), and that $f_0$ was proved to be generic by Herrero and Vel\'azquez in \cite{HVasns92} and \cite{HVcras92} (the proof was given only in one space dimension, and the authors asserted the proof holds also in higher dimensions).

When $\mu\neq 0$ and $q<q_c$, Ebde and Zaag \cite{EZ} were able to show that the existence of a solution with the same profile $f_0$ as for the case $\mu=0$. This is reasonable, since in similarity variables defined below by \eqref{framwork:selfsimilar}, the gradient term comes with an exponentially decreasing term. However, some involved parabolic regularity arguments were needed in \cite{EZ} to handle the gradient term.

\smallskip
When $\mu\neq 0$ and $q=q_c$, up to our knowledge, there is only one result proving the existence of blow-up
solutions for equation \eqref{eq:uequation}: if $\mu<0$ and $p-1$ is small, Souplet, Tayachi and Weissler constructed a selfsimilar blow-up solution in \cite{STWesbugiumj96}. Let us also mention the numerical result by Nguyen in \cite{Y} who finds the same behavior as in \eqref{profile0} with almost the same profile as $f_0$, in the sens that only the constant $b_0$ changes into $b_\mu$, continuous in terms of $\mu$ (let us also mention the solution by Galaktionov and V\'azquez in \cite{GV1} and \cite{GV2} in the supercritical case $q=2>q_c$ with $\mu>0$).

\smallskip
In \cite{TZ15} where we consider the critical case $q=q_c$,  we initially wanted to prove rigourously the numerical result by Nguyen, but we didn't succeed. We ended instead by finding a {\it new} type of blow-up behavior, different from \eqref{profile0}, in the case \eqref{condition:h}, as we state in the following result:

\begin{Thmk}{(\normalfont Blow-up profile for Equation \eqref{eq:uequation})}
\label{th1}
 For any $\varepsilon>0,$ Equation \eqref{eq:uequation} has a
solution $u(x,t)$ such that $u$ and $\nabla u$ blow up in finite
time $T>0$ simultaneously at the origin and only there. Moreover:
\begin{itemize}
\item[(i)] For all $t\in [0,T),$
\begin{eqnarray}
&&\left\|(T-t)^{{\frac 1{p-1}}}u(x,t)-\left(p-1+\frac{b|x|^2}{(T-t)|\log(T-t)|^\beta}\right)^{-{\frac 1{p-1}}}
\right\|_{L^{\infty}(\R^N)}\label{bupprof1}\\ &&\hspace{2cm}\leq {C\over 1+|\log(T-t)|^{\min\left({2\over p-1},{p-3\over 2(p-1)}\right)-\varepsilon}},\nonumber
\end{eqnarray}
 where
 \begin{equation}
\label{bbupprof1}
\beta={p+1\over 2(p-1)},\; \quad \; 
b={1\over 2}(p-1)^{{p-2\over p-1}}\left((4\pi)^{N\over 2}(p+1)^2N\over p
\int_{\R^N}|y|^qe^{-{|y|^2/4}}dy\right)^{{p+1\over p-1}}\mu^{-{p+1\over p-1}}>0,\;
\end{equation}
and a similar estimate holds for $\nabla u$.
 \item[(ii)] For all $x\not=0,$ $u(x,t)\to u^*(x)$ as $t\to T$ in $C^1\left(\frac 1R <|x|<R\right)$ for any $R>0$,
 with
\begin{equation*}
u^*(x)\sim \left({b|x|^2\over \left[2\left|\log|x|\right|\right]^{{p+1\over p-1}}}\right)^{-{1\over p-1}},\; \mbox{ as } \; x\to\; 0,
\end{equation*}
and for $|x|$ small,
 $|\nabla u^*(x)|\leq  C {|x|^{-{p+1\over
p-1}}\over |\log|x||^{\alpha}},$
 for some $\alpha=\alpha(p,\epsilon)\in \R$.
 \end{itemize}
\end{Thmk}
 
\begin{rem}{\rm
Note that the  solution constructed in the above theorem does not exist in the case of the standard nonlinear heat equation, i.e. when $\mu=0$ in \eqref{eq:uequation}. Indeed, our solution has a profile depending on the reduced
variable
$$z={x\over \sqrt{T-t}\left|\log(T-t)\right|^\beta}$$ whereas, we know from the results in \cite{HViph93,Vtams93} that the blow-up profiles in the case $\mu=0$ depend on the reduced variables
$$z={x\over \sqrt{T-t}\left|\log(T-t)\right|^{{1\over 2}}}\; \mbox{or}\; z={x\over (T-t)^{{1\over 2m}}},\; \mbox{where}\; m\geq 2\; \mbox{is an integer}.$$}
\end{rem}

As a consequence of our techniques, we also obtained the following stability result in \cite{TZ15}:

\begin{Thmk}{} \label{th2} The constructed solution is stable with respect to initial data. 
\end{Thmk}
 
Let us give an idea of the methods used to prove the results. We construct the blow-up solution with the profile in Theorem \ref{th1}, by following the methods of  \cite{BK} and \cite{MZsbupdmj97}, though we are far from a simple adaptation, since the gradient term needs genuine new ideas as we explain shortly below.
This kind of methods has been applied for various nonlinear evolution equations.
 For hyperbolic equations, it has been successfully used for the construction of multi-solitons for the semilinear wave equation in one space dimension (see \cite{CZ}). For parabolic equations, it has been used in \cite{MZ08} and \cite{Zihp98}
 for the complex Ginzburg-Landau equation with no gradient structure. See also the cases of the wave maps in \cite{RD}, the Schr\"odinguer maps in \cite{MRR}, the critical harmonic heat follow in \cite{RS1}, the two-dimensional Keller-Segel equation in \cite{RS2} and the nonlinear heat equation involving a subcritical nonlinear gradient term in \cite{EZ}.
 Recently, this method has been applied for a non variational parabolic system in \cite{NZ} and for a logarithmically perturbed nonlinear heat equation in \cite{NguyenZ}.

\smallskip
 Unlike in the subcritical case in \cite{EZ}, the gradient term in the critical case induces substantial changes in the blow-up profile as we pointed-out in the comments following Theorem \ref{th1}. Accordingly, its control requires special arguments.
 So, working in the framework
 of \cite{MZsbupdmj97}, some crucial modifications are needed. In particular, we have to overcome the following challenges:
 \begin{itemize}
  \item[-] The prescribed profile is not known and not  obvious to find. See Section \ref{secform} for a formal approach to justify such a profile, and the introduction of the parameter $\beta$ given by \eqref{equbeta} below.
  \item[-] The profile is different from the profile in \cite{MZsbupdmj97}, hence also from all the previous studies in the parabolic case (\cite{MZsbupdmj97,EZ,NguyenZ,NZ}).
Therefore, brand new estimates are needed. See Section \ref{secexis} below.
   \item[-] In order to handle the new parameter $\beta$ in the profile, we introduce a new shirking set to trap the solution. See Definition \ref{prop:V-def} below. Finding such a set is not trivial,
  in particular the limitation $p>3$ in related to the choice of such a  set.
   \item[-] A good understanding of the dynamics of the linearized operator of equation \eqref{eq:wequation} below around the new profile is needed, taking into account the new shrinking set.
  \item[-] Some crucial global and pointwise estimates of the gradient of the solution as well as fine parabolic regularity results are needed. 
 \end{itemize}
 Then, following \cite{MZsbupdmj97}, the proof is divided in two steps. First, we reduce the problem to a finite dimensional one. Second, we solve the finite dimensional problem and conclude by contradiction, using index theory.

\smallskip
The stability result, Theorem \ref{th2}, is proved similarly as in
\cite{MZsbupdmj97} by interpreting the finite dimensional problem in
terms of the blow-up time and the blow-up point.

 \smallskip
Thanks to simple change of variables, we obtain similar statements for the following perturbation of the following viscous Hamilton-Jacobi (vHJ) equations: 
\begin{equation}\label{eqvHJ}
\partial_t u  = \Delta u +|\nabla u|^q+\nu|u|^{p-1}u,\;
\mbox{ with } \nu>0, \; 3/2<q<2,\; p={q\over 2-q}.
\end{equation}
\begin{COROA}{\normalfont (Blow-up in the viscours Hamilton-Jacobi (vHJ) equation)}Theorems \ref{th1} and \ref{th2} yield stable blow-up solutions in equation \eqref{eqvHJ}.
Moreover, 
the solution and its gradient blow up simultaneously and only at one point. The blow-up profile is given by (\ref{bupprof1}) with appropriate scaling.
\end{COROA}
{
{\begin{rem}
{\rm  Obviously, our result does not hold for the viscous Hamilton-Jacobi equation with $\nu=0$. An interesting question is to understand the behavior of the constructed solutions, say $u_\nu$, as $\nu \to \infty$. In our opinion, this is a difficult open question}
\end{rem}
}}

This note is ogranized as follows:\\
- In Section 2, we explain formally how we obtain the profile and the exponent $\beta$;\\
- In Section 3, we give a formulation of the problem in order to justify the formal argument;\\
- In Section 4, we give the proof of the existence of the profile assuming some technical results.\\
For simplicity, we only focus on the case 
\[
N=1,
\]
and refer the reader to \cite{TZ15} where the high-dimensional case is presented better. We also refer to \cite{TZ15} for the technical details which are omitted here.

\section{A Formal Approach}\label{secform}
\setcounter{equation}{0}
The aim of this section is to explain formally how we derive the behavior given in Theorem \ref{th1}. In particular, how we obtain the profile 
$\vpz$ in \eqref{bupprof1}  (see \eqref{feq} for the notation $\vpz$), the parameter $b$ and the exponent $\beta=2(p+1)/(p-1)$ in \eqref{bbupprof1}.
{
We will also explain why our strategy works only for $\mu>0$, as asserted in \eqref{condition:h}, and not in the case $\mu<0$ (of course, we never consider the case $\mu=0$ which corresponds to the well-known semilinear heat equation). For that purpose, we only assume here that
\[
\mu\neq 0,
\]
and we will explain at the end of this section why we need the positivity assumption on $\mu$ (see \eqref{mu>0} below).
}

\bigskip
 
Let us consider an  arbitrary $T>0$ and the self-similar transformation of
(\ref{eq:uequation})
\begin{equation}\label{framwork:selfsimilar}
w(y,s)= (T-t)^{\frac{1}{p-1}}u(x,t),\; y=\frac{x}{\sqrt{T-t}},~
  s=-\log{(T-t)}.
  \end{equation}
It follows that if $u(x,t)$ satisfies (\ref{eq:uequation}) for all
$(x,t)\in \R\times [0,T),$ then $w (y,s)$ satisfies  the following
equation:
\begin{equation}\label{eq:wequation}
 \partial_s w=\partial_y^2w-\frac{1}{2}y\partial_y
  w-\dfrac{1}{p-1}w+\mu|\partial_y
  w|^{q}+|w|^{p-1}w,
\end{equation}
for all $(y,s)\in \R\times[-\log T,\infty).$ Thus, constructing a
solution $u(x,t)$ for the equation (\ref{eq:uequation}) that blows
up at $T<\infty$ like $(T-t)^{-{\frac{1}{p-1}}}$ reduces to
constructing  a global solution $w(y,s)$ for equation
(\ref{eq:wequation}) such that
\begin{equation}
\label{225}
 0<\varepsilon\leq
\limsup_{s\to\infty}\|w(s)\|_{L^\infty(\R)}\leq {1\over \varepsilon}.
\end{equation}

A first idea to construct a blow-up solution for \eqref{eq:uequation}, would be to find a stationary solution for (\ref{eq:wequation}), yielding a self-similar solution for \eqref{eq:uequation}.
It happens that when $\mu<0$ and $p$ is close to $1,$ the first author together with Souplet and Weissler were able in \cite{STWesbugiumj96}   to construct such a solution. Now, if $\mu>0,$ we know, still from \cite{STWesbugiumj96}  that it is not possible to construct
such a solution in some restrictive class of solutions (see \cite[Remark 2.1, p. 666]{STWesbugiumj96}), of course, apart from the trivial constant solution $w\equiv \kappa$ of (\ref{eq:wequation}), where
\begin{equation} \label{equkappa}\kappa=\Big({1\over p-1}\Big)^{{1\over p-1}}.
\end{equation}

\subsection{Inner expansion}

Following the approach of Bricmont and Kupiainen in \cite{BK}, we may look for a solution $w$ such that $w\to \kappa$ as $s\to \infty.$ Writing $$w=\kappa+\overline{w},$$ we see that
$\overline{w}\to 0$ as $s\to \infty$ and   satisfies the equation:
\begin{equation}
 \label{wbarre}
\partial_s\overline{w} =\mathcal{L}\overline{w}+\overline{B}(\overline{w})+\mu|\nabla \overline{w}|^q,
\end{equation}
 where
\begin{equation}\label{eq:Operator:L}\mathcal{L}=\partial_y^2-\frac{1}{2}y\partial_y+1,~\end{equation}
and
\begin{equation}\label{eqbgare01}
 \overline{B}(\overline{w})= |\overline{w}+\kappa|^{p-1}(\overline{w}+\kappa)-\kappa^p-p\kappa^{p-1}\overline{w}.
\end{equation}
Note that
$$|\overline{B}(\overline{w})-{p\over 2\kappa}\overline{w}^2|\leq C |\overline{w}^3|,$$
where $C$ is a positive constant.

Let us recall some properties of $\mathcal{L}$. The operator $\mathcal{L}$  is  self-adjoint in $D(\mathcal L)
\subset L_\rho^2(\mathbb R)$ where
$$L_\rho^2(\mathbb{R})=\left\{f\in L_{loc}^2(\mathbb R)\; \Big|\; \int_{\mathbb R} \left(f(y)\right)^2\rho(y)dy<\infty \right\}$$
and
$$\rho(y)
=\frac{e^{\frac{-|y|^2}{4}}}{\sqrt{4\pi}},\; y\in \R.$$
 The spectrum of
$\mathcal{L}$ is explicitly given by
$$spec(\mathcal L)=\left\{1-\frac{m}{2}\; \Big| \; m\in \mathbb N\right\}.$$
It consists only in  eigenvalues,
which are all simple,  
and the eigenfunctions
are dilations of Hermite polynomials: the eigenvalue $
1-\frac{m}{2}$ corresponds  to  the following eigenfunction:
\begin{equation} \label{h_m} h_m(y)=\sum_{n=0}^{[\frac{m}{2}]} \frac{m!}{n!(m-2n)!}(-1)^n
y^{m-2n}.\end{equation} In particular $h_0(y)=1,\; h_1(y)=y$ and
$h_2(y)=y^2-2.$ Notice that $h_m$ satisfies:
\[
\int_{\R} h_n h_m \rho dx=2^nn!\delta_{nm}
\mbox{ and }
\mathcal{L}h_m=\Big(1-{m\over 2}\Big)h_m .
\]
   In compliance with the spectral properties of $\mathcal{L},$ we may look for a solution expanded as follows:
\[
\overline{w}(y,s)=\sum_{m\in\N}\overline{w}_m (s) h_m(y).
 \]
 Since $h_m,$ for $m\geq 3$
correspond to negative eigenvalues of $\mathcal{L}$, assuming $\overline{w}$  even in $y,$
we may consider that
\begin{equation}\label{ansatz}
 \overline{w}(y,s)=\overline{w}_0(s)+\overline{w}_2(s)h_2(y),
\end{equation}
with $\overline{w}_0,\; \overline{w}_2\;\to 0$ as $s\to \infty.$

Projecting Equation \eqref{wbarre}, and writing $\mu |\nabla \overline{w}|^q=\mu 2^q|y|^q|\overline{w}_2|^q,$ we derive the following ODE system for $\overline{w}_0$ and $\overline{w}_2:$
$$\overline{w}_0'=\overline{w}_0+{p\over 2\kappa}\left(\overline{w}_0^2+8\overline{w}_2^2\right)+\tilde{c_0}|\overline{w}_2|^q+O\left(|\overline{w}_0|^3+|\overline{w}_2|^3\right),$$
$$\overline{w}_2'=0+{p\over \kappa}\left(\overline{w}_0\overline{w}_2+4\overline{w}_2^2\right)+\tilde{c_2}|\overline{w}_2|^q+O\left(|\overline{w}_0|^3+|\overline{w}_2|^3\right),$$
where
\begin{equation*}
 {\tilde{c}}_0=\mu2^q\int_{\R}|y|^q\rho\; \mbox{ and } \; {\tilde{c}}_2=\frac{\mu2^q}8\int_{\R}|y|^q(|y|^2-2)\rho.
\end{equation*}
Note that for this calculation, we need to know the values of
\begin{equation*}
\int_{\R}(|y|^2-2)^2\rho(y) dy=
 8
\mbox{ and }
\int_{\R}(|y|^2-2)^3\rho(y) dy =
 64.
 \end{equation*}
Note also that the sign of ${\tilde{c}}_0$ and  ${\tilde{c}}_2$ is the same as for $\mu.$ Indeed, obviously $\int_{\R^N}|y|^q\rho(y)dy>0,$ and for $\int_{\R}|y|^q(|y|^2-2)\rho(y)dy,$ using 
integration by parts, we write
\begin{align}
\nonumber
\frac{8\tilde c_2}{2^q\mu}=&\int_{\R}|y|^q(|y|^2-2)\rho(y)dy
  =
\int_{\R}|y|^{q+2}\rho(y)dy-2\int_{\R}|y|^q\rho(y)dy\\ 
 =&
2(q+1)\int_{\R}|y|^{q}\rho(y)dy-2\int_{\R}|y|^q\rho(y)dy \label{8+losange}
 = 2q\int_{\R}|y|^{q}\rho(y)dy>0.
\end{align}
From the equation on $\overline{w}_2',$ we write
$$\overline{w}_2'=\tilde{c_2}|\overline{w}_2|^q\left(1+O\left(|\overline{w}_2|^{2-q}\right)\right)+{p\over \kappa}\overline{w}_0\overline{w}_2+O\left(|\overline{w}_0|^3\right),$$
and assuming that
\begin{equation}
 \label{Hypotesesprofileh1}
|\overline{w}_0\overline{w}_2|\ll |\overline{w}_2|^q,\; |\overline{w}_0|^3\ll |\overline{w}_2|^q,
\end{equation}
we get that
$$\overline{w}_2'\sim sign(\mu)|{\tilde{c}}_2||\overline{w}_2|^q,$$
with $sign(\mu)=1$ if $\mu>0$ and $-1$ if $\mu<0.$\\
In particular, if $\mu>0,$ then $\overline{w}_2$ is increasing tending to $0$ as $s\to \infty$ hence $\overline{w}_2<0$,  while if $\mu<0,$  $\overline{w}_2$ is decreasing tending to $0$ as $s\to \infty$, hence $\overline{w}_2>0.$
Then, since $1<q<2,$
we get
$$\overline{w}_2 \sim -sign(\mu){B\over s^{{1\over q-1}}},$$
with
\begin{equation}
 B  = \left[(q-1)|\tilde{c}_2|\right]^{-{1\over q-1}} \label{5++++triangle}  =  \left[2^{q-2}q(q-1)|\mu|\int_\R |y|^q\rho\right]^{-{1\over q-1}}
\end{equation}
from
\eqref{8+losange}.

From the equation on $\overline{w}_0',$ we write
$$\overline{w}_0'=\overline{w}_0\left(1+O\left(\overline{w}_0\right)\right)+\tilde{c_0}|\overline{w}_2|^q\left(1+O\left(|\overline{w}_2|^{2-q}\right)\right),$$
and assuming that
\begin{equation}
 \label{Hypotesesprofileh2}
|\overline{w}_0'|\ll \overline{w}_0,\;  |\overline{w}_0'|\ll |\overline{w}_2|^q,
\end{equation}
we derive that
$$\overline{w}_0 \sim -{\tilde{c}}_0|\overline{w}_2|^q \sim  {-{\tilde{c}}_0B^q\over s^{{q\over q-1}}}\ll |\overline{w}_2|.$$
Such $\overline{w}_0$ and $\overline{w}_2$ are compatible with the hypotheses \eqref{Hypotesesprofileh1} and \eqref{Hypotesesprofileh2}.

Therefore, since $w=\kappa+\overline{w}$,
it follows from \eqref{ansatz} that
\begin{eqnarray}
\nonumber
w(y,s) & = &\kappa+ \overline{w}_2(s)(|y|^2-2)+o\left(\overline{w}_2\right)\\ \nonumber
&=&\kappa-\frac{sign(\mu)}{s^{\frac 1{q-1}}}B(|y|^2-2)+o\left({1\over s^{{1\over q-1}}}\right)\\ \label{9++losange}
&=&\kappa
-sign(\mu)B{|y|^2\over s^{\frac 1{q-1}}}+2\frac{sign(\mu)}{s^{\frac 1{q-1}}}B+o\left(\frac 1{s^{\frac 1{q-1}}}\right),
\end{eqnarray}
in $L^2_\rho(\R),$ and also uniformly on compact sets by standard parabolic regularity.

\subsection{Outer expansion} From \eqref{9++losange}, we see that the variable
$$z={y\over s^\beta},\; \mbox{ with } \beta={1\over 2(q-1)}={p+1\over 2(p-1)},$$
as given in \eqref{bbupprof1}, is perhaps the relevant variable for
blow-up. Unfortunately, \eqref{9++losange} provides no shape, since
it is valid only on compact sets (note that $z\to 0$ as $s\to
\infty$ in this case). In order to see  some shape, we may need to
go further in space, to the ``outer region'', namely when $z\not=
0.$ In view of \eqref{9++losange}, we may try to find an expression
of $w$ of the form
\begin{equation}
 \label{1++++triangle}
 w(y,s)= \vpz(z)+{a\over s^{2\beta}}+O\left({1\over s^\nu}\right),
\end{equation}
for some $\nu>2\beta.$ Plugging this ansatz in equation \eqref{eq:wequation}, keeping only the main order, we end-up with the following equation on $\vpz:$
\begin{equation}
 \label{ordre0profil}
-{1\over 2}z[\vpz]'(z)-{1\over
p-1}\vpz(z_0)+[{\vpz}(z)]^p=0,\; z={y\over
s^\beta}.\end{equation} Recalling that our aim is to find $w$ a
solution of \eqref{eq:wequation} such that $w\to \kappa$ as $s\to
\infty$ (in $ L^2_\rho$, hence uniformly on every compact set), we
derive from \eqref{1++++triangle} (with $y=z=0$) the natural
condition
$$\vpz(0)=\kappa.$$
Recalling also that we already adopted radial symmetry for the inner equation, we do the same here. Therefore, integrating equation \eqref{ordre0profil}, we see that
\begin{equation}
 \label{2++++triangle}
 \vpz(z)=\Big(p-1+b|z|^2\Big)^{-{1\over
p-1}},
\end{equation}
for some $b\in \R.$ Recalling also that we want a solution $w\in
L^\infty(\R),$ (see \eqref{225}),  we see that $b\geq 0$ and for a
nontrivial solution, we should have
\begin{equation}
 \label{4++++triangle}
 b>0.
\end{equation}
Thus, we have just obtained from \eqref{1++++triangle} that
\begin{equation}
 \label{3++++triangle}
 w(y,s) =
 \Big(p-1+b|z|^2\Big)^{-{1\over
p-1}}+\frac a{s^{2\beta}}+O\left(\frac 1{s^\nu}\right),\; \mbox{ with }\; z={y\over s^\beta}\mbox{ and }\nu>2\beta.
\end{equation}
We should understand this expansion to be valid at least on compact
sets in $z,$ that is for $|y|< R s^\beta,$ for any $R>0.$

\subsection{Matching asymptotics} Since \eqref{3++++triangle} holds for $|y|< R s^\beta,$ for any $R>0$, it holds also uniformly on compact sets, leading to the following expansion for $y$ bounded:
\begin{equation*}
 w(y,s)
 = \kappa-{\kappa b \over (p-1)^2}{|y|^2\over s^{2\beta}}+{a\over s^{2\beta}}+O\left({1\over s^\nu}\right).
\end{equation*}
Comparing with \eqref{9++losange}, we find the following values for $b$ and $a:$
$$b=sign(\mu){B(p-1)^2\over \kappa } \; \mbox{ and }\; a=2 sign(\mu)B.$$
In particular, from \eqref{4++++triangle} we see that
\begin{equation}\label{mu>0}
\mu>0.
\end{equation}

In conclusion, using  \eqref{5++++triangle}, we see that we have just derived the following profile for $w(y,s):$
$$w(y,s) \sim \varphi (y,s)$$
with
\begin{equation}\label{feq}\varphi(y,s)=\vpz\Big({y\over s^\beta}\Big)+{a\over
s^{2\beta}}:=\Big(p-1+b{|y|^2\over s^{2\beta}}\Big)^{-{1\over
p-1}}+{a\over s^{2\beta}}.
\end{equation}
\begin{equation} \label{equbeta}
\beta={p+1\over 2(p-1) },
\end{equation}
\begin{equation}\label{aeq} a={2b\kappa\over (p-1)^2},\end{equation}
\begin{equation} \label{equb}b={1\over 2}(p-1)^{{p-2\over p-1}}\left(2\sqrt \pi(p+1)^2\over p
\int_{\R}|y|^qe^{-{|y|^2/4}}dy\right)^{{p+1\over p-1}}\mu^{-(p+1)/(p-1)},
\end{equation}

\section{Formulation of the problem}
\setcounter{equation}{0}
In this section we formulate the problem in order to justify the formal approach given in the previous section. 
  Let $w,\; y$ and $s$ be as in (\ref{framwork:selfsimilar}). Let us
introduce $\vq(y,s)$ such that
\begin{equation}\label{eq:q+phi+def} w(y,s)=\varphi(y,s)+\vq(y,s),\end{equation}
where $\varphi$ is given by \eqref{feq}.
If $w$ satisfies the equation (\ref{eq:wequation}), then $\vq$
satisfies the following equation:
\begin{equation}\label{qequation}
\partial_s \vq=(\mathcal{L}+V)\vq+B(\vq)+G(\vq)+{R}(y,s),
\end{equation}
where $\mathcal{L}$ is defined by \eqref{eq:Operator:L} and
\begin{equation}\label{eq:V}
 V(y,s)=p~\varphi(y,s)^{p-1}
-\frac{p}{p-1},
\end{equation}
\begin{equation}\label{eq:B}
B(\vq)=|\varphi+\vq|^{p-1}(\varphi+\vq)-\varphi^p-p\varphi^{p-1}\vq,\end{equation}
\begin{equation}\label{RNeq}
\begin{split}
&R(y,s)=\partial_y^2\varphi-\frac{1}{2}y\partial_y\varphi-\frac{\varphi}{p-1}+
\varphi^p-\frac{\partial \varphi}{\partial s}+\mu |\partial_y \varphi|^q
\end{split}\end{equation}
and
 \begin{equation}\label{G} G(v)=\mu|\partial_y\varphi+\partial_y \vq|^q-\mu|\partial_y\varphi|^q.
 \end{equation}
 Our aim is to construct  initial data $v(s_0)$ such that the equation
(\ref{qequation}) has a solution $\vq(y,s)$ defined for all
$(y,s)\in \R\times[-\log T,\infty),$ and satisfies:
\begin{equation}\label{eq:qt0} \lim_{s\rightarrow \infty
}\|\vq(s)\|_{W^{1,\infty}(\R)}=0.\end{equation} From Equation
(\ref{feq}), one sees that the variable  $z={y\over s^\beta}$ plays
a fundamental role. Thus we will consider the dynamics for $|z|>K$
and $|z|<2K$ separately for some $K>0$ to be fixed large.
Since
 \begin{equation}\label{estim:B-R}
|B(\vq)|\leq C|\vq|^2,\;  \|R(.,s)\|_{L^\infty}\leq \f{C}{s}, \|G(\vq)\|_{L^\infty(\R)}\leq {C\over \sqrt{s}} \|\vq\|_{L^\infty(\R)},\end{equation}  for $s$ large enough 
(see \cite{TZ15}), 
it is then reasonable to think that the
dynamics of equation (\ref{qequation}) are influenced by the linear
part, namely $\mathcal{L}+V$.

\bigskip
The properties of the operator $\mathcal{L}$ were given in Section 2. In particular, $\mathcal{L}$ is predominant on all the modes, except on the null modes where the terms $V\vq $ and $G(\vq)$ will play a crucial role 
(see \cite{TZ15}).

As for the potential  $V,$ it has two
fundamental properties  which will strongly  influence our strategy:
\begin{itemize}
 \item[(i)] we have  $ V(.,s)\rightarrow 0$ in
$L_\rho^2(\mathbb{R})$ when $s\rightarrow \infty.$ In practice,  the
effect of $V$ in the blow-up area $(|y|\leq Cs^\beta)$ is regarded
as a perturbation of the effect of  $\mathcal{L}$ (except on the null mode).
 \item[(ii)] outside of the blow-up area, we have the following property:
for all $\epsilon >0$, there exists $C_{\epsilon}>0$ and
$s_{\epsilon}$ such that $$ \sup_{s\geq
s_{\epsilon},~\frac{|y|}{s^{\beta}}\geq
C_{\epsilon}}\left|V(y,s)-(-\frac{p}{p-1})\right|\leq \epsilon,$$
with $-\frac{p}{p-1}<-1.$ As $1$ is the largest eigenvalue of the
operator $\mathcal{L}$, outside the blow-up area  we can consider
that  the operator $ \mathcal {L}+V$ is  an operator with negative
eigenvalues, hence, easily controlled.
\end{itemize}

Considering the fact that the behavior of $V$ is not the same inside
and outside  the blow-up area, we decompose  $ \vq$ as follows. Let
us consider a non-increasing cut-off function  $\chi_0\in
C^\infty\big([0,\infty),[0,1]\big)$ such that
$\mbox{supp}(\chi_0)\subset[0,2]$ and $\chi_0 \equiv 1 $ in $[0,1]$,
and introduce
\begin{equation}\label{def:chi}
\chi(y,s)=\chi_0\left(\frac{|y|}{K~s^{\beta}}\right)\end{equation}
with 
$K$ is some large enough constant so that various estimates in the proof hold.\\
Then, we write
\begin{equation}\label{def:q:proj}
\vq(y,s)=\vq_b(y,s)+\vq_e(y,s),
\end{equation}
with
\begin{equation}\label{def:q:projbis}
\vq_b(y,s)=\vq(y,s)\chi(y,s) \mbox{ and } \vq_e(y,s)=
\vq(y,s)\big(1-\chi(y,s)\big).\end{equation} We remark that $$
\mbox{supp }\vq_b(s)\subset B(0,2K{s^{\beta}}),\; \mbox{supp}~
\vq_e(s)\subset \mathbb{R}\setminus B(0,Ks^{\beta}).$$
As for $\vq_b,$ we will decompose it according to the sign of the eigenvalues of $\mathcal{L}$, by writing
\begin{equation}
\label{def:q_b}
\vq_b(y,s)= \sum_{m=0}^2
\vq_m(s)h_{m}(y)+\vq_{-}(y,s),
\end{equation}
where for $0\le m\le 2$, $v_m=P_m(v_b)$ and $\vq_{-}(y,s)=
P_{-}(\vq_{b})$, with $P_m$ the $L^2_\rho$ projector on $h_m$, the eigenfunction corresponding to $\lambda=1-\frac m2\ge 0$, and $P_{-}$  the projector on $\{h_i,\;|\;i\ge 3\}$, the negative
subspace of the  operator $\mathcal{L}$
 (as announced in the beginning of the section, hereafter, we assume that $N=1$ for simplicity).\\
 Thus, we can decompose $\vq$ in  five components  as follows:
\begin{equation}\label{qprojection}
\vq(y,s)=\sum_{m=0}^2
\vq_m(s)h_{m}(y)+\vq_{-}(y,s)+\vq_e(y,s).
\end{equation}
Here and
throughout the paper, we call $\vq_{-} $ the negative mode of $\vq$,
$\vq_2 $ the null mode of $\vq$, and the subspace spanned by
$\left\{h_m\; |\; m\geq 3\right\}$ will be referred to as the negative
subspace.

\section{The existence proof without technical details}\label{secexis}
\setcounter{equation}{0}
In this section, we prove the existence of a solution $\vq$ of  (\ref{qequation}) such that
\begin{equation}
 \label{but1}
\lim_{s\to \infty}\|\vq(s)\|_{W^{1,\infty}(\R)}=0.
\end{equation}
This is in fact the main step towards the proof of Theorem \ref{th1}. Here,
 we only give the arguments of the proof, and for the technical details, we refer the interested reader to our paper \cite{TZ15}. For the remaining steps of the proof of Theorem \ref{th1} and also for the proof of Theorem \ref{th2}, we refer to \cite{TZ15}.

\smallskip
Since $p>3,$ we see that, by definition of $\beta$ given by \eqref{equbeta}, $\beta\in ({1\over 2},1).$
Our construction is build on a careful choice of the initial data for $\vq$ at a time $s_0.$ We will choose it in the following form:

\begin{defi}{\normalfont Choice of the initial data)}\label{initial data}
Let us define, for $A\geq 1,$ $s_0=-\log T>1$ and $d_0,\;d_1\in \R,$
the function
\begin{equation}\label{initial-data}
 \psi_{s_0,d_0,d_1}(y)={A\over s_0^{2\beta+1}}\Big(d_0h_0(y)+d_1h_1(y)\Big)\chi(2y,s_0),
\end{equation}
where $h_i,\; i=0,\,1$ are defined by (\ref{h_m}) and $\chi$ is
defined by (\ref{def:chi}).
\end{defi}

The solution of equation \eqref{qequation} will be denoted by
$\vq_{s_0,d_0,d_1}$ or $\vq$ when there is no ambiguity. We will
show that if $A$ is fixed large enough, then, $s_0$ is fixed large
enough depending on $A,$  we can fix the parameters $(d_0,d_1)\in
[-2,2]^2,$ so that the solution $\vq_{s_0,d_0,d_1}(s)\to 0$ as $s\to
\infty$ in $W^{1,\infty}(\R),$ that is, \eqref{but1} holds. Owing to
the decomposition given in \eqref{initial-data}, it is enough to
control the solution in a shrinking set  defined as follows:
\begin{defi}{\normalfont  (A set shrinking to zero)}\label{prop:V-def} Let $\gamma$ be any real number  such that
\begin{equation}\label{gamma}
3\beta<\gamma<\min(5\beta-1,2\beta+1).
\end{equation}
 For all  $ A\geq 1$ and $
s\geq 1$, we define $\vartheta_A(s)$ as the set of all functions $r
\in L^\infty(\R)$ such that

\[||r_e||_{L^\infty(\R)}~\leq~ {A^2\over
s^{\gamma-3\beta}},\; \;\;\;\Big\|{r_-(y)\over 1+|y|^3}
\Big\|_{L^\infty(\R)}~\leq~ {A\over s^{\gamma}},\; ~\]

\[|r_0|,\; |r_1|\leq
~{A\over s^{2\beta+1}},\;\;\;~|r_2|~\leq~{\sqrt{A}\over
s^{4\beta-1}} ,\]
\\
where $r_-$, $r_e$ and  $r_m$ are defined in (\ref{qprojection}).
\end{defi}
\begin{rem}
{\rm Since $p>3,$ it follows that $\frac 12<\beta<1,$ in particular the range for $\gamma$ in \eqref{gamma} is not empty. Of course, the set $\vartheta_A(s)$ depends also on the choice of $\gamma$ satisfying \eqref{gamma}.
However, while $A$ will be chosen large enough so that various estimates hold, $\gamma$
will be fixed once for all throughout the proof.}
\end{rem}

Since $A\geq1,$  then the sets $\vartheta_A(s)$ are increasing (for
fixed $s$) with respect to $A$ in the sense of inclusion.  We also
show  the following property of elements of $\vartheta_A(s):$

 \smallskip

{\it For all $A\geq 1$, there exists $s_{01}(A) \geq 1$ such that, for
all $s\geq s_{01}$  and $r\in \vartheta_A(s),$ we have
\begin{equation}
 \label{prop:V-defbis}
 \|r\|_{L^\infty(\R)}\leq C\frac{A^2}{s^{\gamma-3\beta}},
\end{equation}
where  $C$ is a positive constant 
(see \cite{TZ15}).
 }

\smallskip
By \eqref{prop:V-defbis}, if a solution $v$ stays in $\vartheta_A(s)$ for $s\geq s_0,$ then it converges to $0$ in $L^{\infty}(\R)$ (the convergence of the gradient will follow from parabolic regularity). Reasonably, our aim is then reduced to prove the following proposition:

\begin{Prop}{\normalfont (Existence of solutions trapped in  $\vartheta_A(s)$)}
\label{dans vs} There exists $A_2\ge 1$ such that for $A\ge A_2$ there exists $s_{02}(A)$ such that for all $s_0\ge s_{02}(A)$, there exists  $(d_0,d_1)$ such that if
$v$ is the solution of (\ref{qprojection}) with initial data at $s_0,$  
given by  (\ref{initial-data}), then
 $\vq(s)\in \vartheta_A(s)$, for all $s\geq s_0.$
\end{Prop}

This proposition gives the stronger convergence to $0$ in
$L^{\infty}(\R)$ thanks to \eqref{prop:V-defbis}, and the
convergence in $W^{1,{\infty}}(\R)$ will follow from an involved parabolic regularity argument as explained in \cite{TZ15}.

Let us first make sure that we can choose the initial data such that it starts in  $\vartheta_A(s_0).$ In other words, we will define a set where we will at the end select the good parameter $(d_0, d_1)$ that will give the conclusion of Proposition \ref{dans vs}. More precisely,
we have the following result:

\begin{Prop}{\normalfont (Properties of initial data)}\label{prop:diddc} For each $A\geq 1$, there
exists $s_{03}(A)>1$ such that for all $s_0\geq s_{03}$,
there exists a rectangle
\begin{equation*}
\mathcal{D}_{s_0}\subset [-2,2]^2
 \end{equation*}
 such that the mapping
\begin{eqnarray*}
\R^2&\rightarrow & \R^2,\\(d_0,d_1) &\mapsto&
(\psi_0,\psi_1).
\end{eqnarray*}
(where $\psi$ stands for $\psi_{s_0,d_0,d_1}$) is linear, one to
one from $\mathcal{D}_{s_0}$ onto $[-{A\over s_0^{2\beta+1}},{A\over
s_0^{2\beta+1}}]^2$ and maps $
\partial \mathcal{D}_{s_0} $ into $\partial\Big(
[-{A\over s_0^{2\beta+1}},{A\over s_0^{2\beta+1}}]^2\Big)$. Moreover, it has degree one on the boundary.
 
\end{Prop}

\smallskip
\noindent
{\bf Proof.}  See \cite{TZ15}.

\smallskip
\noindent
{\bfseries Proof of Proposition \ref{dans vs}}. Let us consider $A\geq 1,$ $s_0\geq s_{03},$ $(d_0,d_1)\in \mathcal{D}_{s_0},$ where $s_{03}$ is given by Proposition \ref{prop:diddc}. From the existence theory (which follows from the Cauchy problem
for equation \eqref{eq:uequation} in $W^{1,\infty}(\R))$ mentioned in the introduction), starting in $\vartheta_A(s_0)$ which is in $\vartheta_{A+1}(s_0),$ the solution stays in $\vartheta_A(s)$ until some maximal time $s_*=s_*(d_0,d_1).$ If
 $s_*(d_0,d_1)=\infty$ for some $(d_0,d_1)\in \mathcal{D}_{s_0}$, then  the proof is complete. Otherwise, we argue by contradiction and suppose that $s_*(d_0,d_1)<\infty$ for any $(d_0,d_1)\in \mathcal{D}_{s_0}.$
 By continuity and the definition of $s_*$,  the solution at the  point $s_*$, is on the boundary of
 $\vartheta_{A}(s_*).$ Then, by definition of $\vartheta_{A}(s_*),$ one at least of the inequalities in that definition  is an equality. Owing to the following proposition, this can happen only
 for the first two components. Precisely, we have the following result:

 \begin{Prop}{\normalfont (Control of  $\vq(s)$ by $(\vq_0(s),\vq_1(s))$ in $ \vartheta_A(s)$)}\label{prop:rt2dimen} There
exists $A_4\geq 1$ such that  for each $A\geq~A_4,$
there exists
$s_{04}(A)\in \R$
such that for all $s_0\geq s_{04}(A)$,
 the following holds:

If $\vq$ is a solution of (\ref{qequation}) with initial data at $s=s_0$ given by
(\ref{initial-data}) with $(d_0,d_1)~\in~\mathcal{D}_{s_0} $, and $\vq(s)~\in~\vartheta_A(s)$  for all  $ s \in~[s_0,s_1],~$
with $\vq(s_1)\in \partial \vartheta_A(s_1)$ for some $s_1\geq s_0$, then:
\begin{itemize}
\item[(i)] (Reduction to a finite dimensional problem) We have: $$ \left(\vq_0(s_1),\vq_1(s_1)\right)\in \partial\left( \left[-\f{A}{s_1^{2\beta+1}},\f{A}{s_1^{2\beta+1}}\right]^2\right).$$
\item[(ii)] (Transverse crossing) There exist $m\in \{0,1\}$ and $\omega \in \{-1,1\}$ such that
\begin{equation*}
 \omega \vq_m(s_1)=\f{A}{s_1^{2\beta+1}}\mbox{ and }  \omega\vq_m'(s_1)>0.
\end{equation*}
\end{itemize}
\end{Prop}

\noindent
Assume the result of the previous proposition, for which the proof is given in \cite{TZ15}, and continue the proof of
Proposition \ref{dans vs}. Let $A\ge A_4$ and $s_0\geq s_{04}(A).$
It follows from Proposition \ref{prop:rt2dimen},
part (i), that $\left(\vq_0(s_*),\vq_1(s_*)\right)\in \partial\left( \left[-\f{A}{s_*^{2\beta+1}},\f{A}{s_*^{2\beta+1}}\right]^2\right),$
and the following function \begin{align*}\Phi&: \mathcal{D}_{s_0} \rightarrow \partial\left([-1,1]^2\right)\\
&(d_0,d_1)\mapsto
\dfrac{s_*^{2\beta+1}}{A}(\vq_0,\vq_1)_{(d_0,d_1)}(s_*),\; \mbox{ with } s_*=s_*(d_0,d_1),\end{align*}
is well defined. Then, it follows from Proposition \ref{prop:rt2dimen}, part (ii) that  $\Phi$ is continuous. On the other hand, using Proposition \ref{prop:diddc}, parts (i) and (ii)  together with  the fact that
$\vq(s_0)= \psi_{s_0,d_0,d_1}$, we see that when $(d_0,d_1)$ is on the boundary of the rectangle $\mathcal{D}_{s_0},$ we have strict
inequalities for the other components. Applying the transverse crossing property  given in Proposition \ref{prop:rt2dimen}, part (ii), we see that
$\vq(s)$  leaves  $\vartheta_A(s)$ at $s=s_0$,  hence $s_*(d_0,d_1)=s_0$. Using Proposition \ref{prop:diddc}, part (i), we see that the  restriction of $\Phi$
to the boundary is of degree 1. A contradiction then follows from the index theory. Thus, there exists a value $(d_0,d_1)\in \mathcal{D}_{s_0}$
such that for all $s\geq s_0,~ \vq_{s_0,d_0, d_1}(s)\in \vartheta_A(s)$.
This concludes the proof of  Proposition \ref{dans vs}.

\bigskip
{\it Completion of the proof of \eqref{but1}.} By Proposition
\ref{dans vs} and \eqref{prop:V-defbis}, it remains only to show that
$\|\nabla \vq(s)\|_{L^\infty(\R)}\to 0$ as $s\to \infty.$ This in fact follows from a very involved parabolic regularity argument 
given in \cite{TZ15}, which implies that there exists $s_{05}$ such that for  $s\geq s_{05}$,
$$\|\vq(s)\|_{W^{1,\infty}(\R)}\leq {C(A)\over s^{\gamma-3\beta}},$$
hence, by \eqref{gamma}, \eqref{but1} follows by taking $s_{02}\geq \max\left(s_{01},s_{03},s_{04},s_{05}\right)$.

\bibliographystyle{plain}

\begin{thebibliography}{99}

\small 
 \bibitem{AW92}
{\sc L.~Alfonsi {\normalfont and} F.~B. Weissler}, {Blow up in {${\bf R}\sp n$} for a parabolic equation with a damping
  nonlinear gradient term}, In Nonlinear diffusion equations and their equilibrium states, 3
  (Gregynog, 1989), volume~7 of  Progr. Nonlinear Differential Equations
  Appl., pages 1-20. Birkh\"auser Boston, Boston, MA, 1992.

\bibitem{BergerKhon}
{\sc   M. Berger  {\normalfont and}  R. V. Koh}, {A rescaling algorithm for the numerical calculation of blowing-up solutions}, Comm. Pure Appl. Math., 41 (1988), 841-863

\bibitem{BK}
{\sc J.~Bricmont  {\normalfont and}  A.~Kupiainen}, {Universality in blow-up for nonlinear heat equations}, Nonlinearity, 7 (1994), 539-575.

\bibitem{CWbngsj89}
{\sc  M.~Chipot  {\normalfont and}  F.~B. Weissler} , {Some blowup results for a nonlinear parabolic equation with a
  gradient term}, SIAM J. Math. Anal., 20 (1989), 886-907.

\bibitem{CFQ}
{\sc  M.  Chlebik,  M.  Fila  {\normalfont and}  P.  Quittner} , {Blowup of positive
 solutions of  a semilinear parabolic equation with a gradient term}, Dyn. Contin.
Discrete Impulsive Syst. Ser. A, {\bf 10} (2003), pp. 525-537.

\bibitem{CZ}
{\sc R. C\^ote  {\normalfont and}  H. Zaag}, {Construction of a multisoliton blowup solution to the semilinear wave equation in one space dimension}, Comm. Pure Appl. Math., 66 (2013), 1541-1581

\bibitem{EZ}
{\sc M. A. Ebde  {\normalfont and}  H. Zaag}, {Construction and stability of a blow up solution for a nonlinear heat equation with a gradient term}, SEMA J., 55 (2011), 5-21

 
\bibitem{GP1}
{\sc V. A. Galaktionov,  {\normalfont and}  S. A. Posashkov}, {The equation {$u_t=u_{xx}+u^\beta$}. Localization, asymptotic behavior of unbounded solutions}, Akad. Nauk SSSR Inst. Prikl. Mat. Preprint, 97 (1985), 30 pages.

\bibitem{GP2}
{\sc V. A. Galaktionov,  {\normalfont and}  S. A. Posashkov}, {Application of new comparison theorems to the investigation of unbounded solutions of nonlinear parabolic equations},
Differentsialcnye Uravneniya, 22 (1986),1165-1173, 1285. English translation:  Differential Equations, 22 (1986), 809-815.

\bibitem{GV1}
{\sc V. A. Galaktionov  {\normalfont and}  J. L. Vazquez}, {Regional blow-up in a semilinear heat equation with convergence to a Hamilton-Jacobi equation}, SIAM J. Math. Anal., 24 (1993), 1254-1276.

\bibitem{GV2}
{\sc V. A. Galaktionov  {\normalfont and}  J. L. Vazquez}, {Blow-up for quasilinear heat equations described by means of nonlinear Hamilton-Jacobi equations}, J. Differential Equations, 127 (1996), 1-40.
 

\bibitem{HVcras92}
{\sc  M.~A. Herrero  {\normalfont and}  J.~J.~L. Vel{\'a}zquez}.
\newblock Comportement g\'en\'erique au voisinage d'un point d'explosion pour
  des solutions d'\'equations paraboliques unidimensionnelles.
\newblock {C. R. Acad. Sci. Paris S\'er. I Math.}, 314(3):201-203, 1992.

\bibitem{HVasns92}
{\sc M.~A. Herrero  {\normalfont and}  J.~J.~L. Vel{\'a}zquez}.
\newblock Generic behaviour of one-dimensional blow up patterns.
\newblock {Ann. Scuola Norm. Sup. Pisa Cl. Sci. (4)}, 19(3):381-450, 1992.

\bibitem{HViph93}
{\sc M. A. Herrero   {\normalfont and} d J. J. L. Vel{\'a}zquez}, {Blow-up behavior of one-dimensional
semilinear parabolic equations}, Ann. Inst. Henri Poincar\'e, 10 (1993),  131-189.

 
\bibitem{MZ08}
{\sc N.~Masmoudi  {\normalfont and}  H.~Zaag}, {Blow-up profile for the complex {G}inzburg-{L}andau equation}, J. Funct. Anal., 255 (2008), 1613-1666.

 
\bibitem{MRR}
{\sc F. Merle, P. Rapha\"el  {\normalfont and}   I. Rodnianski}, {Blow up dynamics for smooth equivariant solutions to the energy critical {S}chr\"odinger map}, C. R. Math. Acad. Sci. Paris, 349 (2011), 279-283.

\bibitem{MZsbupdmj97}
{\sc  F.~Merle  {\normalfont and}   H.~Zaag}, {Stability of the blow-up profile for equations of the type {$u\sb  t=\Delta u+\vert u\vert \sp {p-1}u$}}, Duke Math. J., 86 (1997), 143-195.

\bibitem{NguyenZ}
{\sc V. T. Nguyen  {\normalfont and}   H. Zaag}, {Construction of a stable blow-up solution for a class of strongly perturbed semilinear heat equations}, arXiv:1406.5233, (2014). submitted.

\bibitem{NZ}
{\sc N. Nouaili  {\normalfont and}    H. Zaag}, {Profile for a simultaneously blowing up solution for a complex valued semilinear heat equation}, Comm. Partial Differential Equations,
40 (2015), 1197-1217. 


\bibitem{RD}
{\sc P. Rapha{\"e}l,   {\normalfont and}  I. Rodnianski}, {Stable blow up dynamics for the critical co-rotational wave  maps and equivariant {Y}ang-{M}ills problems}, Publ. Math. Inst. Hautes \'Etudes Sci., (2012), 1-122.

\bibitem{RS1}
{\sc P. Rapha{\"e}l  {\normalfont and}   R. Schweyer}, {Stable blowup dynamics for the 1-corotational energy critical harmonic heat flow}, Comm. Pure Appl. Math., 66 (2013), 414-480.

\bibitem{RS2}
{\sc  P. Rapha{\"e}l  {\normalfont and}   R. Schweyer}, {On the stability of critical chemotactic aggregation}, Math. Ann.,359 (2014),267-377.

 
 \bibitem{SnT}
 {\sc S.  Snoussi  {\normalfont and}  S.  Tayachi},  {Large time behavior of solutions for parabolic equations with nonlinear gradient terms.} Hokkaido Mathematical Journal, 36 (2007), 311-344.

\bibitem{SnTW}
{\sc S.  Snoussi, S.  Tayachi  {\normalfont and}  F.~B. Weissler},  {Asymptotically self-similar global solutions of a semilinear parabolic equation with a nonlinear gradient term.} Proc. Roy. Soc. Edinburgh Sect. A, 129 (1999),
1291-1307.

\bibitem{Smmas96}
{\sc P. Souplet}, {Finite time blow-up for a nonlinear parabolic equation with a gradient term and applications},
Math. Methods in the Applied Sciences, 19 (1996), 1317-1333.

\bibitem{Sejde2001}
{\sc P. Souplet}, {Recent results and open problems on parabolic equations with gradient nonlinearities},
Electron. J. Differential Equations 10 (2001), 19 pp.

\bibitem{Snepp2005}
{\sc P. Souplet},  {The influence of gradient perturbations on blow-up asymptotics in semilinear parabolic problems: a survey},
Nonlinear elliptic and parabolic problems, 473-495, Progr. Nonlinear Differential Equations Appl., 64, Birkh\"auser, Basel, 2005.

\bibitem{ST}
{\sc P. {Souplet}  {\normalfont and}  S. {Tayachi}}, {Blowup  rates for
nonlinear heat equations with gradient terms and for parabolic inequalities},  Colloq. Math, 88 (2001), 135-154.


\bibitem{STWesbugiumj96}
{\sc  P.~Souplet, S.~Tayachi,  {\normalfont and}  F.~B. Weissler}, {Exact self-similar blow-up of solutions of a semilinear parabolic
  equation with a nonlinear gradient term}, Indiana Univ. Math. J., 45 (1996), 655-682.

\bibitem{SW}
{\sc P. Souplet  {\normalfont and}   F. B.  Weissler}, {Self-Similar Subsolutions and Blowup for Nonlinear Parabolic Equations}, Journal of Mathematical Analysis and Applications,  212 (1997),  60-74.

\bibitem{SW1}
{\sc P. Souplet  {\normalfont and}   F. B.  Weissler}, {Poincar\'e
 inequality and global solutions of  a
nonlinear parabolic equation}, Ann. Inst. Henri Poincar\'e, Analyse
non lin\'eaire,  16 (1999),  337-373.

\bibitem{TZ15}
{\sc S. Tayachi  {\normalfont and}  H. Zaag}, {Existence of a stable blow-up profile for the nonlinear heat equation with a critical power nonlinear gradient term}, (2015), arXiv:1506.08306, submitted.


\bibitem{Vtams93}
{\sc J. J. L. Vel{\'a}zquez}, {Classification of singularities for blowing up solutions in higher dimensions},
Trans. Amer. Math. Soc., 338 (1993), 441-464.

 
\bibitem{Y}
{\sc V. T. Nguyen}, {Numerical analysis of the rescaling method for parabolic problems with blow-up in finite time}, arXiv:1403.7547,  (2014). submitted

\bibitem{Zihp98}
{\sc  H.~Zaag}, {Blow-up results for vector-valued nonlinear heat equations with no
  gradient structure}, Ann. Inst. H. Poincar\'e Anal. Non Lin\'eaire, 15(1998), 581-622.

\end {thebibliography}

\end{document}